\documentstyle[12pt]{article}
\title{The time constant  vanishes only on the percolation cone  in  directed   first passage percolation
\footnotetext{AMS classification: 60K 35.}
\footnotetext{Key words and phrases: directed first passage percolation, growth model, and phase transition.}} 
\author{Yu Zhang\footnote{Research supported by NSF grant DMS-4540247. }\\University of Colorado}
\date{}
\begin{document}
\baselineskip .20in
\maketitle

\begin{abstract}
We consider the directed first passage percolation model on ${\bf Z}^2$.
In this model, we assign  independently to each edge $e$ a passage time $t(e)$
with a common distribution $F$.
We denote by
$\vec{T}({\bf 0}, (r,\theta))$  
the passage time from the origin to $(r, \theta)$ by a northeast path for 
$(r, \theta)\in {\bf R}^+\times [0,\pi/2]$. It is known that 
$\vec{T}({\bf 0}, (r, \theta))/r$ converges to a time constant $\vec{\mu}_F (\theta)$.
Let $\vec{p}_c$ denote the critical probability for oriented percolation.
In this paper, we show that the time constant has a phase transition  divided by $\vec{p}_c$, as follows:

 (1) If  $F(0) < \vec{p}_c$, then $\vec{\mu}_F(\theta) >0$ for all
$0\leq \theta\leq \pi/2$. 

(2) If  $F(0) = \vec{p}_c$, then $\vec{\mu}_F(\theta) >0$ if and only if $\theta\neq \pi/4$. 

(3) If $F(0)=p > \vec{p}_c$, then there exists a percolation cone between $\theta_p^-$ and $\theta_p^+$ for
 $0\leq \theta^-_p< \theta^+_p \leq  \pi/2$ such that $\vec{\mu} (\theta) >0$ if and only if $\theta\not\in [\theta_p^-, \theta^+_p]$. Furthermore, all the moments of $\vec{T}({\bf 0}, (r, \theta))$ converge whenever $\theta\in [\theta_p^-, \theta^+_p]$.\\
As applications, we describe the 
shape  of the directed growth model on the distribution of $F$. We  give a phase transition for the shape
divided by $\vec{p}_c$.

\end{abstract}
\section {Introduction of the model and results.}
In this  directed first passage percolation model, 
we consider the vertices of  the ${\bf Z}^2$ lattice  and the  edges of the vertices with the Euclidean distance 1.
We assign independently to each edge a non-negative {\em passage time} $t(e)$ with a common distribution $F$.
More formally, we consider the following probability space. As the
sample space, we take $\Omega=\prod_{e\in {\bf Z}^2} [0,\infty),$ whose points 
are called {\em configurations}.
Let ${\bf P}=\prod_{e\in {\bf Z}^2} {\mu}_e$ be the corresponding product measure on $\Omega$, where
${\mu}_e$ is the measure on $[0, \infty)$ with distribution $F$. The
expectation  with respect to ${\bf P}$ is denoted by ${\bf E}(\cdot)$.
For any two vertices $u$ and $v$ in ${\bf Z}^2$, 
a  path $\gamma$ from $u$ to $v$ is an alternating sequence 
$(v_0, e_1, v_1,...,v_i, e_{i+1}, v_{i+1},...,v_{n-1},e_n, v_n)$ of vertices $v_i$ and 
  edges $e_i$ between $v_i$ and $v_{i+1}$ in ${\bf Z}^2$, with $v_0=u$ and $ v_n=v$. 
For a vertex $u$, its northeast edges from $u$ are denoted by 
$u=(u_1, u_2)$ to $(u_1+1, u_2)$ or  to $(u_1, u_2+1)$.
Given a path $(v_0, e_1, v_1,...,v_i, e_{i+1}, v_{i+1},...,v_{n-1},e_n, v_n)$,  
if each edge $e_i$ is a northeast edge from $v_i$, the path
is called  {\em northeast}, or {\em directed}.
For short, we denote northeast edges or northeast paths by NE edges or NE paths.

Given  a path $\gamma$, we define its passage time  as 
$$T(\gamma)=\sum_{e\in \gamma} t(e).$$
For any two vertices $u$ and $v$, we define the passage time from  $u$ to $v$ by
$${T}(u, v)=\inf \{ {T}(\gamma)\},$$
where the infimum is over all possible  paths from  $u$ to $v$. 
We also define 
$$\vec{T}(u,v)=\inf \{ {T}(\gamma)\},$$
where the infimum is over all possible NE paths from  $u$ to $v$. 
If there does not exist
a NE path from $u$ and $v$, we simply define
$$\vec{T}(u,v)=\infty.$$
A NE path $\gamma$ from $u$ to $v$ with ${T}(\gamma)=\vec{T} (u,v)$
 is called an {\em optimal path} of   $\vec{T} (u,v)$. 
We need to point out that the optimal path may  not be  unique.
If we focus on a special configuration $\omega$, we may write $\vec{T} (u,v)(\omega)$,
instead of $\vec{T} (u,v)$.
 
In addition to vertices on ${\bf Z}^2$, we may also consider points on ${\bf R}^2$.
In particular, 
we often use the polar coordinates $\{(r, \theta)\}= {\bf R}^+\times [0, \pi/2]$, 
where $r$ and $\theta$ represent the radius 
and the  angle between the radius and the $X$-axis, respectively.
We may extend the definition of passage time over ${\bf R}^+\times [0, \pi/2]$.
If $u=(r, \theta)$ in ${\bf R}^+\times [0, \pi/2]$, 
we define $\vec{T} ({\bf 0}, u)=\vec{T} ({\bf 0}, u')$, where
$u'$  is the nearest neighbor of $u$  in
${\bf Z}^2$. Possible indetermination can be eliminated by choosing an
order on the vertices of ${\bf Z}^2$ and taking the smallest nearest
neighbor for this order. Similarly,   $\vec{T}(u, v)$ can be defined for
any $u, v\in {\bf R}^2$. Moreover, with this extension, for any points $u$ and $v$ in ${\bf R}^2$,  we may consider a path 
of ${\bf Z}^2$ from $u$ to $v$.

Given a non-zero vector $(r, \theta)\in {\bf R}^+\times [0, 2\pi]$, by a subadditive argument, 
if $Et(e) < \infty$, then
$$
\lim_{r\rightarrow \infty}{1\over r} {T}({\bf 0}, (r, \theta)) = \lim_{r\rightarrow \infty}{1\over r} {\bf E}{T}({\bf 0}, (r, \theta))=
\inf_{r} {1\over r} {\bf E} {T}({\bf 0}, (r, \theta))=\mu_F( \theta) \mbox{ a.s. and in } L_1.\eqno{}
$$
We call $\mu_F (\theta)$ a  {\em time constant}.
Furthermore, by the same subadditive argument, for a non-zero vector $(r, \theta)\in {\bf R}^+\times [0, \pi/2]$,
$$
\lim_{r\rightarrow \infty}{1\over r} \vec{T}({\bf 0}, (r, \theta)) = \lim_{r\rightarrow \infty}{1\over r} {\bf E}\vec{T}({\bf 0}, (r, \theta))=
\inf_{r} {1\over r} {\bf E} \vec{T}({\bf 0}, (r, \theta))=\vec{\mu}_F( \theta) \mbox{ a.s. and in } L_1.\eqno{(1.1)}
$$
We also call $\vec{\mu}_F (\theta)$   a  time constant. By the subadditive argument again, we know 
(see Proposition 2.1 (iv) in Martin (2004)) that
$$\vec{\mu}_F (\theta)\mbox{  is finite and convex in }\theta.\eqno{(1.2)}$$
 In general, we require that $t(e)$ has a finite first moment or  $m$-th moment. However,
we sometimes require  the following  stronger tail assumption:
$${\bf E}\exp(\eta \tau(e))< \infty  \mbox{ for } \eta >0.\eqno{(1.3)}$$

Recall the undirected first passage percolation model for $\{T(u, v)\}$. Kesten (1986) showed  that 
there is a phase transition divided by critical probability, $p_c$, of bond percolation for a time constant. More precisely, he showed that  time constant $\mu_F(\theta)$ vanishes  if and only if $F(0) \geq p_c$. Therefore, $F(0) > p_c$, $F(0)=p_c$, and
$F(0) < p_c$ are called  the {\em supercritical},
the {\em critical},  and the {\em subcritical} phases, respectively.
It is natural to examine a similar  situation
for the directed first  passage percolation model. In this paper, our focus is  that there
is also a  phase transition for $\vec{\mu}_F(\theta)$ divided by
critical probability, $\vec{p}_c$, of directed bond percolation.  
We will demonstrate for the supercritical and critical phases, which are quite different
from the undirected first passage percolation model (see Kesten and Zhang (1997), and Zhang (1995)).
We will also examine the subcritical phase, which is similar to the undirected model (see Kesten (1986)).\\

{\bf 1.1. Supercritical phase.} We now focus  on the supercritical phase: $F(0) > \vec{p}_c.$
Before introducing our results, we would like to introduce a few basic oriented
percolation results.
If we rotate our lattice counterclockwise by $45^\circ$ and  extend each edge by a factor of $\sqrt{2}$, 
the new graph is 
denoted by ${\cal L}$ with  
oriented edges from $(m,n)$ to $(m+1, n+1)$ and
to $(m-1, n+1)$. 
Each edge is independently open or closed with probability $p$ or $1-p$. 
An oriented path from $u$ to $v$ is defined as a sequence $v_0=u, v_1,\cdots, v_m=v$ of points of ${\cal L}$. The path 
has the vertices $v_i=(x_i,y_i)$ and $v_{i+1}=(x_{i+1}, y_{i+1})$ for $0\leq i\leq m-1$ such that $y_{i+1} =y_i+1$ and $v_i$
and $v_{i+1}$ are connected by an oriented  edge. An oriented path is open if each of its edges is open.
For two vertices $u$ and $v$ in ${\cal L}$,  we say $u\rightarrow v$ if there is an oriented open path
from $u$ to $v$. For $A\subset (-\infty, \infty)$,
we denote a random subset by
$$\xi_n^A=\{x: \, \exists \,\,\, x'\in A \mbox{ such that } (x',0) \rightarrow (x,n) \}\mbox{ for } n>0.$$
The  right edge for this set is defined by
$$r_n=\sup \xi_n^{(-\infty, 0]} \,\,\, (\sup \emptyset =-\infty).$$
By a subadditive argument (see section 3 (7) in Durrett (1984)),
there exists a non-random constant $\alpha_p$ such that
$$ \lim_{n\rightarrow \infty} {r_n\over n}=\lim_{n}\frac{{\bf E} r_n}{ n}=\alpha_p \mbox{ a.s. and in }L_1,
\eqno{(1.4)}$$
where $\alpha_p >0$ if $p> \vec{p}_c$, and $\alpha_p=0$ if $p=\vec{p}_c$, and $\alpha_p=-\infty$ if $p< \vec{p}_c$.
Now we rotate the lattice back to ${\bf Z}^2$. If $p\geq \vec{p}_c$, the  {\em percolation cone} is 
the cone between
  two polar equations $\theta=\theta_p^\mp$ in the first quadrant, where (see Marchand (2002))
$$\theta_p^\mp=\arctan \left( { 1/2\mp \alpha_p/\sqrt{2} \over  1/2\pm \alpha_p/\sqrt{2} }\right ).$$
Note that if $p=\vec{p}_c$, then the percolation cone  shrinks to the positive diagonal line. 
In fact, for any point $(r, \theta)$ with $\theta\in [\theta_p^-,   \theta_p^+]$, it can be shown 
(see Lemma 3 in Yukich and Zhang (2006)) that
$${\bf P}[\exists \mbox{ a NE zero-path from the origin to $(r, \theta)$}]> C.\eqno{(1.5)}$$
In this paper,  $C$ and $C_i$ are always positive  constants that may depend on $F$, but not on $t$, $r$, $k$, or $n$. 
Their values 
are not significant and  may change from 
appearance to appearance.
With these definitions, we have the following theorem regarding the passage time on the  percolation
cone:\\

{\bf Theorem 1.} {\em If $F(0)=p >\vec{p}_c $ and ${\bf E}(t(e))^m < \infty$ for $m \geq 1$, then 
for all $r$ and $\theta\in [\theta_p^- , \theta_p^+]$, there exists $C=C(F,m)$ such that}
$${\bf E}\vec{T}({\bf 0}, (r, \theta))^m\leq C.$$

In contrast to the passage time on the  percolation cone, we have another theorem:\\

{\bf Theorem 2.} {\em If $F(0)=p >\vec{p}_c$ and $\theta\not\in [\theta_p^-,  \theta_p^+]$, then 
for all $r$, there exist $\delta=\delta(F, \theta)>0$ and $C_i=C_i(F, \theta, \delta)$ for $i=1,2$ such that}
$${\bf P}[ \vec{T}({\bf 0}, (r, \theta)) \leq \delta r]\leq C_1 \exp(-C_2 r).$$

Together with Theorems 1 and 2, we have the following corollary:\\

{\bf Corollary 3.} {\em If $F(0)=p >\vec{p}_c$ and ${\bf E}(t(e)) < \infty$, then for $0\leq \theta \leq \pi/2$,}
$$\vec{\mu}_F(\theta)=0 \mbox{ iff } \theta \in [\theta_p^- , \theta_p^+].$$

{\bf Remark 1.} We would like to discuss  $\vec{\mu}_F(\theta)$ as a function of $F$.
 Recall that in the general first passage percolation model, Yukich and Zhang (2006) showed that the time constant
is not third differentiable in the direction of $\theta_p^\pm$. We find out that the same proof together with
$$T(u, v) \leq \vec{T}(u, v)$$
can be carried out to show the same result
for directed first passage percolation. 
Here we state the following result but omit the proof.
We denote by $\vec{\mu}_F(\theta, p)$ the time constant for $F(0)=p$.
If $t(e)$ only takes two values 0 or 1 and $F(0) >\vec{p}_c$, then
$$\mbox{ $\vec{\mu}_F(\theta_p^\pm, p)$ is not
third differentiable in $p$.}\eqno{(1.6)}$$
Except in these two directions, we believe that there is no other singularity.\\


{\bf Remark 2.} Note that $\vec{\mu}_F(\theta)$ can also be considered as a function of $\theta$. By the convexity in (1.2), 
we can show that $\vec{\mu}_F(\theta)$ is continuous in $\theta$. 
We believe that $\theta_p^\mp $ are also the singularities
for $\vec{\mu}_F(\theta)$ in $\theta$.\\

{\bf Conjecture 1.} If $F(0) > \vec{p}_c$, show that $\vec{\mu}_F(\theta)$ has  singularities at $\theta^\pm_p$.\\

{\bf 1.2. Critical phase.} We focus on the critical phase:
$F(0) = \vec{p}_c.$  Now, as we mentioned, the percolation cone shrinks to the positive diagonal line.
Similar to the supercritical phase, we can show the following theorem:\\

{\bf Theorem 4.} {\em If $F(0) = \vec{p}_c$  and 
$\theta\neq \pi/4$, then there exist $\delta=\delta(F, \theta)>0$ and $C_i=C_i(F, \theta, \delta)$ for $i=1,2$ 
such that}
$${\bf P}[ \vec{T} ({\bf 0}, (r, \theta))\leq \delta r]\leq C_1 \exp(-C_2 r).$$

The time constant at $\theta=\pi/4$  has  double behaviors: supercritical and subcritical behaviors.  
First, we show that it has a supercritical behavior:\\

{\bf Theorem 5.} {\em If  ${\bf E}t(e) < \infty$ and 
$F(0) = \vec{p}_c$, then 
$$\vec{\mu}_F( \pi/4)=0 .\eqno{(1.7)}$$
In addition, if $F(0)=p$, then}
$$\lim_{p\rightarrow \vec{p}_c} \vec{\mu}_F( \pi/4)=0 .\eqno{(1.8)}$$

{\bf Remark 3.} 
Cox and Kesten used  a circuit method (1981) to show  the following result, which is a stronger result than (1.8). 
If $F_n \Rightarrow F$, 
then 
$$\lim_{n\rightarrow \infty}{\mu}_{F_n}(\theta)={\mu}_{F}(\theta).$$
However, their method cannot be applied for the directed model, since
a path may not be directed after using a piece of circuit. Therefore, we might need a new method to solve
this problem. \\

{\bf Conjecture 2.} If $F_n \Rightarrow F$,  show $\vec{\mu}_{F_n}(\theta)=\vec{\mu}_{F}(\theta)$ for the directed model.\\

Together with Theorems 4 and 5, we have the following corollary:\\

{\bf Corollary 6.} {\em If $F(0)=\vec{p}_c$ and ${\bf E}(t(e)) < \infty$, then for $0\leq \theta \leq \pi/2$,}
$$\vec{\mu}_F(\theta)=0 \mbox{ iff } \theta =\pi/4.$$

To pursue the convergent rate, we need to 
 use the {\em isoperimetric} inequality by Talagrand (1995).   
Denote by ${\cal S}$ the sets of all NE paths from the origin to $(r, \theta)$ with the minimum passage time. Let 
$$\alpha =\sup_{\gamma\in {\cal S}} |\gamma |,$$
where $|A|$ is the number of vertices in $A$ for some vertex set $A$.
Since we only focus on NE paths,
$$\alpha\leq C r.\eqno{}$$
Denote by $M$ a median of $\vec{T}({\bf 0}, (r, \theta))$. By Theorem 8.3.1 (see Talagrand (1995)),
if (1.3) holds, then there exist constants $C_i=C_i(F, \theta)$ for $i=1,2$
such that
$${\bf P}\left[|\vec{T}({\bf 0}, (r, \theta))-M|\geq t\right]\leq C_1\exp\left (-C_2\min\left \{{t^2\over \alpha}, t\right\}\right).\eqno{}$$
By this  isoperimetric inequality together with a simple computation, we can show the following argument.
 For all $r>0$ and $1\leq x\leq \sqrt{r}$, if (1.3) holds, then
$${\bf P}\left[|\vec{T}({\bf 0}, (r, \theta))-{\bf E}\vec{T}({\bf 0}, (r, \theta))|\geq x\sqrt{r}\right]\leq C_1\exp(-C_2 x^2).\eqno{}$$
With this concentration inequality, 
we  can use Alexander's result  (1996)  to show the following. For all $r$, if (1.3) holds, there exists $C= C(F,\theta)$  such that for all $0< r$ 
$$r\vec{\mu}(\theta) \leq {\bf E} \vec{T}({\bf 0}, (r, \theta))\leq r \vec{\mu}_F (\theta)+ C\sqrt{r} \log r.\eqno{(1.9)}$$
With (1.9) and Theorem 5, if (1.3) holds and $F(0)=\vec{p}_c$, then there exists $C=C(F)$ such that
$${\bf E}\vec{T}({\bf 0}, (r, \pi/4))\leq C \sqrt{r} \log r.\eqno{(1.10)}$$

{\bf Remark 4.} The upper bound might  not be tight at the right side of (1.10). In fact, we believe 
the following conjecture in a much tight upper bound:\\

{\bf Conjecture 3.}
If (1.3) holds and $F(0)=\vec{p}_c$,  show that
$${\bf E}\vec{T}({\bf 0}, (r, \pi/4))\leq C \log r.\eqno{(1.11)}$$
Note that (1.11) holds for the undirected first passage time (see Chayes, Chayes, and Durrett (1986)).
In contrast,  the lower bound is more complicated. It might depend on
 how $F(x) \downarrow F(0)=\vec{p}_c$ as $x\downarrow 0$. When the right derivative of $F(0)$ is large enough, 
we believe the following conjecture occurs as  the same as the  undirected model (see Zhang (1995)):\\

{\bf Conjecture 4.} There exists $F$ with $F(0) =p_c$ such that
$$E\vec{T}({\bf 0}, (r, \pi/4))\leq C.\eqno{(1.12)}$$
However, when the right derivative of $F(0)$ is small, we believe that it has a subcritical behavior
 similar to the behavior of undirected passage time. More precisely,
$$\lim_{r\rightarrow \infty} {\bf E}\vec{T}({\bf 0}, (r, \pi/4))=\infty.\eqno{(1.13)}$$
In fact, we may simply  ask the same questions when  $t(e)$ only takes 0 and 1 with $F(0)=\vec{p}_c$.\\

{\bf Conjecture 5.}  If $t(e)$ only takes 0 and 1 with $F(0)=\vec{p}_c$, show that 
$$C_1\log r \leq {\bf E}\vec{T}({\bf 0}, (r, \pi/4))\leq C_2 \log r .\eqno{(1.14)}$$
Note that (1.14) is indeed true (see Chayes, Chayes, and Durrett (1986)) for the undirected critical model.
Furthermore, Kesten and Zhang (1997) showed a central limit theorem for the passage time in
the undirected critical model.
Here, we  partially verify (1.14) for the directed critical model:\\

{\bf Theorem 7.} {\em  If $t(e)$ only takes two values 0 and 1 with $F(0) =\vec{p}_c$, then }
$$\lim_{r\rightarrow \infty}{\bf E}\vec{T}({\bf 0}, (r, \pi/4))=\infty.$$

{\bf Remark 5.} As we mentioned  above, we know  the continuity of $\vec{\mu}_F(\theta)$ in $\theta$. 
We believe that
there is a power law when $\theta\rightarrow \pi/4$.
More precisely, we assume that $F(0) = \vec{p}_c$ and $t(e)$ only takes values 0 and 1.\\

{\bf Conjecture 6.}
$\vec{\mu}_F(\theta)\approx |\theta-\pi/4|^\alpha$
for some $0< \alpha < 1$.\\

{\bf 1.3. Subcritical phase.} Finally, we focus on the subcritical phase:
$F(0) =p < \vec{p}_c.$
On this phase, we show the following theorem:\\

{\bf Theorem 8.} {\em If  $F(0) < \vec{p}_c$, then for all $r$ and $0\leq \theta\leq \pi/2$, there 
exist $\delta=\delta(F)$ and  $C_i=C_i(F, \delta)$ for $i=1,2$  such that}
$${\bf P}[ \vec{T} ({\bf 0}, (r, \theta))\leq \delta r]\leq C_1 \exp(-C_2 r).$$

By Theorem 8, there exists $C=C(F)$ such that for all $r$ and $\theta$
$${\bf E}[\vec{T}({\bf 0}, (r, \theta))] \geq Cr .\eqno{(1.15)}$$
With (1.15) and (1.1), we have the following corollary:\\

{\bf Corollary 9.} {\em If ${\bf E}t(e) < \infty$ and $F(0) < \vec{p}_c$, then for all $0\leq \theta\leq \pi/2$,}
$$\vec{\mu}_F(\theta) >0.$$

{\bf Remark 6.} We would like to focus on a special case in the subcritical phase. 
In fact, Hammersley and Welsh (1965) considered
$t(e)+a$ for some real number $a$. They used $F\oplus a (x)=F(x-a)$ to denote the distribution.
Clearly, if $a >0$, each edge takes at least a time $a$, so $F\oplus a (0) =0$. Therefore, it is in a subcritical phase.
Durrett and Liggett (1981) consider the case that
$$F\oplus a(a) > \vec{p}_c\eqno{(1.16)}$$ 
for undirected passage time $T(u, v)$. If we consider directed passage time $\vec{T}(u, v)$ with (1.16),
note that 
$$ \vec{T}({\bf 0}, (x, y))= \vec{T}'({\bf 0}, (x, y))+ax+ay,$$
where $\vec{T}'(u, v)$ is passage time from $u$ to $v$ with passage time $t(e)$ on edge $e$.
Thus, the directed first passage percolation model on $F\oplus a (x)$  is equivalent to 
the supercritical phase discussed before. \\
 
{\bf 1.4.  Shape of the growth model.}
We may discuss the shape theorem for this directed first passage percolation. Define the shape as
$${\bf C}_t=\{(r, \theta)\in {\bf R}^+\times [0, \pi/2]: \vec{T}({\bf 0}, (r, \theta))\leq t\}.$$
For each $(r, \theta)\in {\bf R}^+\times [0, \pi/2]$, by the subadditive argument,
$$
\lim_{s\rightarrow \infty}{1\over s} \vec{T}({\bf 0}, (s r, \theta)) = \lim_{s\rightarrow \infty}{1\over s}{\bf E} \vec{T}({\bf 0}, (s r, \theta))
=\vec{\mu}_F(r, \theta) \mbox{ a.s. and in } L_1.\eqno{(1.17)}
$$
By (1.1) and (1.17), we know that
$$r\mu_F(\theta)= \mu_F(r, \theta).$$
With (1.17), we define the directed growth shape as
$${\bf C}=\{(r, \theta)\in {\bf R}^+\times [0, \pi/2]: \vec{\mu}_F(r, \theta)\leq 1\}.\eqno{(1.18)}$$
With these definitions,
Martin (2004) proved that if ${\bf E}t^2(e) < \infty$ and 
$$\inf_{r\neq 0} {\vec{\mu}_F(r, \theta) \over r} >0,$$  then ${\bf C}$ is a convex compact set, and for any $\epsilon >0$,
$$(1-\epsilon) {\bf C}\subset {{\bf C}_t\over t} \subset (1+\epsilon ){\bf C}, \mbox{ eventually with probability 1}.\eqno{(1.19)}$$
The result in (1.19) is called  {\em shape theorem}. 
In the subcritical case, for all $0\leq \theta \leq \pi/2$, by Corollary 8, ${\bf C}$ is a convex compact set such that
 the shape theorem 
holds. We denote the shape between two angles by
$${\bf C}(\theta_1, \theta_2) =\{(r, \theta)\in {\bf R}^+\times [\theta_1, \theta_2]: \vec{\mu}_F(r, \theta)\leq 1\},
{\bf C}_t(\theta_1, \theta_2) =\{(r, \theta)\in {\bf R}^+\times [\theta_1, \theta_2]:\vec{T}({\bf 0}, (r, \theta))\leq t \}$$
for $ 0\leq \theta_1 \leq  \theta_2\leq \pi/2$. Furthermore, we denote by 
$\rho(\theta)$ the boundary point $(\rho(\theta), \theta)$ of ${\bf C}$.

In the supercritical case, by Corollary 3 and (1.19),
for any small 
$0< \delta $, 
$${\bf C}(0, \theta_p^--\delta)\mbox{ and } {\bf C}( \theta_p^+ +\delta, \pi/2)\mbox{  are
 convex compact  sets } \eqno{(1.20)}$$
such that the shape theorem holds:
$$(1-\epsilon) {\bf C}(0, \theta_p^--\epsilon)\subset {{\bf C}_t(0, \theta_p^--\delta)\over t} \subset (1+\epsilon ){\bf C}(0, \theta_p^--\delta)$$
and 
$$(1-\epsilon) {\bf C}( \theta_p^++\delta,\pi/2)\subset {{\bf C}_t(\theta_p^++\delta, \pi/2)\over t} \subset (1+\epsilon ){\bf C}( \theta_p^++\delta,\pi/2) \eqno{(1.21)}$$
eventually with probability 1.
On the other hand, by Corollary 3 again,
$${\bf C}(\theta_p^-, \theta_p^+)\mbox{ and } \lim_{t} {{\bf C}_t(\theta_p^-, \theta_p^+)\over t}\mbox{ equal  the percolation cone.} \eqno{(1.22)}$$

In the critical case, 
for any small $0< \delta $, by Corollary 6 and (1.19),
$${\bf C}(0, \pi/4-\delta)\mbox{ and }{\bf C}( \pi/4 +\delta, \pi/2)\mbox{  are
 convex compact sets }\eqno{(1.23)}$$
such that the shape theorem holds:
$$(1-\epsilon) {\bf C}(0, \pi/4-\epsilon)\subset {{\bf C}_t(0, \pi/4-\delta)\over t} \subset (1+\epsilon ){\bf C}(0, \pi/4-\delta)$$
and 
$$(1-\epsilon) {\bf C}( \pi/4+\delta,\pi/2)\subset {{\bf C}_t(\pi/4+\delta, \pi/2)\over t} \subset (1+\epsilon ){\bf C}( \pi/4+\delta,\pi/2) \eqno{(1.24)}$$
eventually with probability 1. On the other hand, by Corollary 6 again, 
$${\bf C}(\pi/4, \pi/4)\mbox{ and }  \lim_{t} {{\bf C}_t(\pi/4, \pi/4)\over t}\mbox{  equal   the positive diagonal line.} \eqno{(1.25)}$$
In particular, in both the supercritical and critical phases,
 by Theorems 1 and 2, and Theorem 5, the continuity of $\mu_F(\theta)$ in $\theta$, 
$$\rho(\theta)\rightarrow \infty \mbox{ as } \theta\rightarrow \theta_p^\pm.\eqno{}$$

\begin{figure}\label{F:graphG}
\begin{center}
\setlength{\unitlength}{0.0125in}%
\begin{picture}(200,200)(-200,800)
\thicklines
\put(0, 800){\vector(1,0){150}}
\put(0, 800){\vector(0,1){150}}
\put(0, 780){$\mbox{Supercritical phase:} F(0) > \vec{p}_c$}

\put(0, 800){\circle*{3}}
\put(0, 800){\line(1,2){50}}
\put(0, 800){\line(2,1){100}}
\put(0, 800){\line(1,3){35}}
\put(0, 800){\line(3,1){100}}

\put(50, 900){$\theta_p^+$}
\put(100, 850){$\theta_p^-$}
\put(20, 920){$\theta_p^++\delta $}
\put(100, 820){$\theta_p^--\delta $}

\put(20, 800){\circle*{4}}
\put(25, 801){\circle*{2}}

\put(30, 802){\circle*{2}}
\put(35, 803){\circle*{2}}
\put(40, 804){\circle*{2}}
\put(45, 805){\circle*{2}}
\put(50, 807){\circle*{2}}
\put(55, 809){\circle*{2}}
\put(60, 811){\circle*{2}}
\put(65, 813){\circle*{2}}
\put(70, 815){\circle*{2}}
\put(75, 817){\circle*{2}}
\put(80, 819){\circle*{2}}
\put(85, 821){\circle*{2}}
\put(88, 825){\circle*{2}}
\put(90, 830){\circle*{3}}
\put(0, 820){\circle*{4}}
\put(1, 825){\circle*{2}}
\put(2, 830){\circle*{2}}
\put(3, 835){\circle*{2}}
\put(4, 840){\circle*{2}}
\put(5, 845){\circle*{2}}
\put(7, 850){\circle*{2}}
\put(9, 855){\circle*{2}}
\put(11, 860){\circle*{2}}
\put(13, 865){\circle*{2}}
\put(15, 870){\circle*{2}}
\put(17, 875){\circle*{2}}
\put(19, 880){\circle*{2}}
\put(21, 885){\circle*{2}}
\put(25, 888){\circle*{2}}
\put(30, 890){\circle*{3}}

\put(-170, 800){\vector(1,0){150}}
\put(-170, 800){\vector(0,1){150}}
\put(-170, 780){$\mbox{Critical phase:} F(0) = \vec{p}_c$}

\put(-170, 800){\circle*{3}}
\put(-170, 800){\line(1,2){50}}
\put(-170, 800){\line(2,1){100}}
\put(-170, 800){\line(1,1){100}}
\put(-70, 900){$\theta=\pi/4$}
\put(-150, 910){$\theta=\pi/4+\delta$}
\put(-90, 855){$\theta=\pi/4-\delta$}
\put(-150, 800){\circle*{4}}
\put(-145, 801){\circle*{2}}
\put(-140, 803){\circle*{2}}
\put(-135, 805){\circle*{2}}
\put(-130, 807){\circle*{2}}
\put(-125, 810){\circle*{2}}
\put(-120, 813){\circle*{2}}
\put(-115, 816){\circle*{2}}
\put(-110, 819){\circle*{2}}
\put(-105, 822){\circle*{2}}
\put(-100, 825){\circle*{2}}
\put(-95, 828){\circle*{2}}
\put(-90, 831){\circle*{2}}
\put(-85, 834){\circle*{2}}
\put(-80, 837){\circle*{2}}
\put(-75, 842){\circle*{2}}
\put(-70, 849){\circle*{4}}

\put(-170, 820){\circle*{4}}
\put(-169, 825){\circle*{2}}
\put(-167, 830){\circle*{2}}
\put(-165, 835){\circle*{2}}
\put(-163, 840){\circle*{2}}
\put(-160, 845){\circle*{2}}
\put(-157, 850){\circle*{2}}
\put(-154, 855){\circle*{2}}
\put(-151, 860){\circle*{2}}
\put(-147, 865){\circle*{2}}
\put(-144, 870){\circle*{2}}
\put(-141, 875){\circle*{2}}
\put(-138, 880){\circle*{2}}
\put(-135, 885){\circle*{2}}
\put(-132, 890){\circle*{2}}
\put(-127, 895){\circle*{2}}
\put(-120, 900){\circle*{4}}

\put(-340, 800){\vector(1,0){150}}
\put(-340, 800){\vector(0,1){150}}
\put(-360, 780){$\mbox{Subcritical phase:} F(0) < \vec{p}_c$}

\put(-340, 820){\circle*{4}}
\put(-339, 825){\circle*{2}}
\put(-338, 830){\circle*{2}}
\put(-336, 835){\circle*{2}}
\put(-333, 840){\circle*{2}}
\put(-330, 845){\circle*{2}}
\put(-328, 850){\circle*{2}}
\put(-326, 855){\circle*{2}}
\put(-324, 860){\circle*{2}}
\put(-320, 865){\circle*{2}}
\put(-315, 870){\circle*{2}}
\put(-310, 875){\circle*{2}}
\put(-305, 880){\circle*{2}}
\put(-292, 885){\circle*{2}}
\put(-284, 887){\circle*{2}}
\put(-276, 889){\circle*{2}}
\put(-267, 890){\circle*{2}}
\put(-259, 887){\circle*{2}}

\put(-320, 800){\circle*{4}}
\put(-315, 801){\circle*{2}}
\put(-310, 802){\circle*{2}}
\put(-305, 804){\circle*{2}}
\put(-300, 807){\circle*{2}}
\put(-295, 810){\circle*{2}}
\put(-290, 813){\circle*{2}}
\put(-285, 815){\circle*{2}}
\put(-280, 820){\circle*{2}}
\put(-275, 825){\circle*{2}}
\put(-270, 830){\circle*{2}}
\put(-265, 838){\circle*{2}}
\put(-260, 846){\circle*{2}}
\put(-258, 852){\circle*{2}}
\put(-256, 860){\circle*{2}}
\put(-254, 869){\circle*{2}}
\put(-253, 876){\circle*{2}}
\put(-259, 887){\circle*{2}}

\end{picture}
\end{center}
\caption{{\em The graph shows  shape ${\bf C}$ in  subcritical, critical, and supercritical phases.
In the supercritical phase, the right figure, the shape is the percolation cone  between two angles $\theta^\pm_p$, and 
the  other two parts of the shape are finite. In the critical phase, the middle figure,
the percolation cone shrinks to the positive diagonal line.
The other two parts of the shape are finite. In the subcritical phase, the left figure, the shape is finite.  }}
\end{figure}
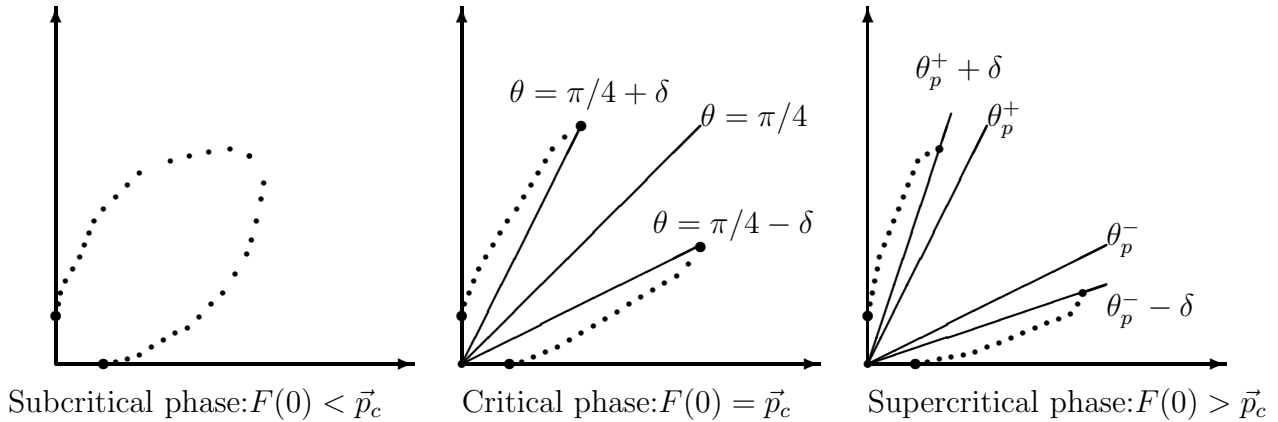 
In the subcritical case, by Corollary 9, the shape theorem in (1.19) holds.  
We can describe the phases of the shapes as  Fig. 1.\\

{\bf Remark 7.} Since the shape is convex, by (1.25),  on the critical and super critical phases, the slope of the line
passing through  $(\rho(\theta_1), \theta_1)$
and $(\rho(\theta_2), \theta_2)$ cannot be more than  $\tan(\theta_p^-)$ when $\theta_1 < \theta_2< \theta_p^-$.
By symmetry, we have the same property when $\pi/2 \geq \theta_1 > \theta_2 > \theta_p^+$.\\


We may relate the directed  first passage percolation to the following directed growth model. At time 1, a cell ${\bf A}_1$  consists of
the unit square with the center at the origin. 
Each square has four edges: the north, the east, the south, and the west 
edges. Two squares are connected if they have a common edge.
Suppose that at time $n$ we have  connected $n$ unit squares, denoted by ${\bf A}_n$. 
Let $\partial {\bf A}_n$ be the boundary of ${\bf A}_n$. A square is a boundary square of ${\bf A}_n$ if one of its edges
belongs to $\partial {\bf A}_n$. We collect all the north and the east edges in $\partial {\bf A}_n$ from the boundary squares. We denote these edges by  the  {\em  northeast
edges} of ${\bf A}_n$. 
At time $n+1$, a new square is added to ${\bf A}_n$ such that it connected to the northeast edges of ${\bf A}_n$.
The location of the new square
is chosen with a probability proportional to the northeast edges of ${\bf A}_n$.

Now we consider $F$ has an exponential distribution with rate 1. Instead of associating the passage time to the edges,
we may associate the passage time to vertices.  By the same discussion, we can define the directed growth shape
${\bf C}_t$ and show the shape theorem of (1.19).
By a similar computation (see page 131 in Kesten (1986)), we can show that  shapes ${\bf A}_n$ and ${\bf C}_{t_n}$
have the same distribution with
$$t_n= \inf\{t: {\bf C}_t \mbox{ contains $n$ vertices}\}.$$
Thus, the shape theorem  for ${\bf C}_t$ implies that $n^{-1/2}{\bf A}_n$ has also an asymptotic shape. 

Unlike the undirected model, the oriented percolation model in higher dimensions has been more limited.
For example, we cannot define the right edge $r_n$ for the oriented percolation model when $d >2$.
However, if one could develop a similar argument of the percolation cone as we defined in section 1.1,
our techniques for first passage percolation would apply for  higher dimensions.
Finally, we conclude this section with the following strictly
monotone conjecture:\\

{\bf Conjecture 7.} If ${\mu}_F(\theta) >0,$ then $\vec{\mu}_F (\theta) > {\mu}_F(\theta)$.

\section{Preliminaries.}

{\bf 2.1. Renormalization method.}
We introduce  the method of renormalization in Kesten and Zhang (1990).
We define, for a large integer $M$ and $w=(w_1, w_2)\in {\bf Z}^2$, the squares by
$$B_M(w)=[Mw_1, Mw_1+M)\times [Mw_2, Mw_2+M) .$$
We denote  these {\em $M$-squares}  by $\{B_M(w): w\in {\bf Z}^2\}$.
 For a path $\gamma$ (not necessary a directed path) starting from the origin, we denote  
a fattened ${\gamma}_M$ by
$${\gamma}_M=\{B_M(w): B_M(w)\cap \gamma \neq \emptyset \}.$$
We denote by $|{\gamma}_M|$ the number of $M$-squares in ${\gamma}_M$.
By our definition, 
$$|\gamma| \geq |\gamma_M| \mbox{ and }  |{\gamma}_M|\geq {|\gamma |\over M^2}.\eqno{(2.1)}$$
For each $M$-square $B_M(w)$, there are eight $M$-square neighbors. We say they are adjacent to $B_M(w)$.
We denote  $B_M(w)$ and its  eight $M$-square neighbors by $\bar{B}_M(w)$.
$\bar{B}M(w)$ is called a {\em $3M$-square}.
Since $\gamma$ is connected, ${\gamma}_M$ has to be connected through the square connections.

If $B_M(w)\cap \gamma\neq \emptyset$ and $\bar{B}_M(w)$ does not contain the origin, note that 
$\gamma$ has to cross $\bar{B}_M(w)\setminus B_M(w)$ to reach to $B_M(w)$, so
$\bar{B}_M(w)$ contains at least $M$ vertices of $\gamma$
in its interior.   We collect all such $3M$-squares $\{ \bar{B}_M(w)\}$ such that
their center $M$-squares contain at least a vertex of $\gamma$. We call these $3M$-squares {\em center $3M$-squares}
of $\gamma$.
With these definitions,  the following lemma (see Zhang, page 22 (2008)) can be calculated directly.
\\

{\bf Lemma 1.} {\em For a connected path $\gamma$, if $|\gamma_M|=k$, then
there are at least $k/15$ disjoint center $3M$-squares of $\gamma$. }\\

{\bf 2.2. Results  for  oriented percolation.}
We assign  either open  or closed to each edge with
probability $p$ or $1-p$ independently from the other edges.
For two sets $A$ and $B$, if there exists a NE open path from $u\in A$ to $v\in B$,
we write  $A\rightarrow B$ as the event.

First, we focus on the subcritical phase: $p < \vec{p}_c$.
Let 
$${\cal C}_0=\{ u: {\bf 0} \rightarrow u\}.$$ 
Durrett (section 7, (6) (1984)) showed the following lemma:\\

{\bf Lemma 2.} {\em If $p < \vec{p}_c$, then there exist $C_i$ for $i=1,2$ such that}
$${\bf P}[ |{\cal C}_0|\geq n ]\leq C_1\exp(-C_2 n).$$

Now we focus on the critical and supercritical phases: $p \geq \vec{p}_c$.
Given two points $u=(u_1, u_2)$ and $v=(v_1, v_2)$,  we define
the  {\em slope} between them by 
$$sl(u, v)={v_2-u_2 \over v_1-u_1}.$$
With these definitions, Zhang (Lemma 3 (2008)) showed  the following lemma:\\

{\bf Lemma 3. } {\em For $0<a< \tan(\theta_p^-)$, 
if $p \geq \vec{p}_c$, then there exist  $C_i=C_i(p,a)$ for $i=1,2$ such that for all $u\in {\bf R}\times [0, \pi/2]$,
$${\bf P}[{\bf 0}\stackrel{}{\rightarrow }u\mbox{ with }sl({\bf 0}, u)\leq a] \leq C_1\exp(-C_2 u_1) . $$}

{\bf 2.3. Analysis for the shape ${\bf C}_t$.}
Now we would like to introduce a few geometric properties for ${\bf C}_t$.
In the remainder of section 2.3, we only consider  $t(e)$ when  it takes value 0 or 1 with $F(0)=p$.
If $t(e)=0$ or 1, $e$ is said to be open or closed.

Given a set $\Gamma\subset {\bf R}^2$, we  let 
$\Gamma'$ be  all vertices on ${\bf Z}^2$ contained in $\Gamma$.
It is easy to see that
$$ \Gamma'\subset \Gamma\subset \{v+(-1, 1)^d: v\in \Gamma'\}.  $$
As we defined in the last section, ${\bf C}_t$  is finite, and so is ${\bf C}'_t$.  
A set  $A$ is said to be 
{\em directly connected} in ${\bf Z}^2$ if any two vertices of $A$ are connected by 
a NE path in $A$.

Given a finite directly connected set $\Gamma$ of ${\bf Z}^2$, 
we  define  its vertex boundary as follows. For each $v\in \Gamma $,
$v\in \Gamma $ is said to be a boundary vertex of $\Gamma$ if there exists $u\not \in \Gamma$ but $u$
is adjacent to $v$ by either a north or an east edge. 
We denote by $\partial \Gamma$  all boundary vertices of $\Gamma$.
 We also let $\partial_o\Gamma $ be  all vertices not in $\Gamma$, but adjacent to $\partial \Gamma$
by north or east edges. $\partial_e \Gamma$ is denoted by these NE edges between $\partial \Gamma$
and $\partial_o \Gamma$.

We define the 
 event as
$$\{{\bf C}'_t=\Gamma\}=\{\omega: {\bf C}'_t(\omega)=\Gamma\}.$$
With these definitions, Zhang (see propositions 1--3 in Zhang (2006))
 proved the following lemmas for undirected first passage percolation.
The proofs can be carried out  by changing paths to directed paths, so we omit the proofs.
In fact, these lemmas are easily understood  by drawing a few figures.\\

{\bf Lemma 4.} {\em ${\bf C}'_t$  is  directly connected.}\\

{\bf Lemma 5.} {\em For all $v\in \partial {\bf C}'_t$, $\vec{T}({\bf 0}, v)=t$, and for all $u\in \partial_o {\bf C}'_t$,
$\vec{T}({\bf 0}, u)=t+1$.}\\

{\bf Lemma 6.} {\em The event of $\{{\bf C}'_t=\Gamma\}$ only depends on the zeros and ones of the edges on 
$\Gamma \cup \partial_o \Gamma $}.\\

{\bf 2.4. Monotone property for the time constant.} 
Finally, we would like to introduce a monotone lemma for the time constant.
Comparing  two distributions $F_1$ and $F_2$, we have the following lemma:\\

{\bf Lemma 7.} {\em If ${\bf E}t(e) < \infty$ and $F_1(x) \leq F_2(x)$ for all $x$, then for $\theta\in [0, \pi/2]$,}
$$ \vec{\mu}_{F_2} (\theta)\leq \vec{\mu}_{F_1}(\theta).$$

{\bf Proof.} Smythe and Wierman (1978) proved the same result in their Theorem 7.12 for undirected 
first passage percolation. The same proof can be carried out to show Lemma 7.  $\Box$\\

\section{ Subcritical phase.} In  section 3, we assume that
$F(0) < \vec{p}_c$. Since $F$ is right-continuous, we take $\epsilon >0$ small such that
$$F(\epsilon)< \vec{p}_c.\eqno{(3.1)}$$
We say that an edge is open if $t(e)\leq \epsilon$, otherwise  $e$ is said to be closed. 
With (3.1), we know that
$${\bf P}[e \mbox{ is open}] < \vec{p}_c.\eqno{(3.2)}$$
Now we work on a NE path from the origin to $(r, \theta)$. As before, we use $\gamma_M$ to denote the squares of $\gamma$.
If a square in $\gamma_M$ contains a closed edge, we call the square  a {\em bad} square.
Otherwise, it is a {\em good} square. If  there is a path $\gamma$ from the origin to $(r, \theta)$
such that it has less than $Cr$ closed edges
for some small $C$, 
then there is less than $C|\gamma_M|$ bad squares in $\gamma_M$. 
Now we account the choices of these squares. Note that $\gamma$ is connected, and so is $\gamma_M$.
 We assume that
$$|\gamma_M|=k.$$
By  Lemma 1 in section 2, we know there are
at least $|\gamma_M|/15$ center $3M$-squares of $\gamma$. Thus, for $C < 1/(30)$, there are at least
$$|\gamma_M|/(15)- C|\gamma_M|\geq |\gamma_M|/(30)\eqno{(3.3)}$$ 
 disjoint 
$3M$-squares such that their center squares contain   an edge of $\gamma$ and all the $M$-squares in these $3M$-squares
 are good. We also call these $3M$-squares  good.
By a standard method (see (4.24) in Grimmett (1999)), there are at most $7^{2k}$ choices for all
possible choices of $\gamma_M$. When $\gamma_M$ is fixed, we select these  good $3M$-squares. There are at most
$2^k$ choices for these good $3M$-squares. 

For each good $3M$-square $\bar{B}_M(w)$,  
there exists a NE open path crossing the $3M$-square from a vertex at the boundary of $B_M(w)$ to another vertex
at the  boundary of $\bar{B}_M(w)$. 
There are at most $4M$ choices for the starting vertex, and the path contains at least $M$ edges. 
For a fixed $\bar{B}_M(w)$,
we denote by ${\cal E}_w$ the event that there exists a NE open path  from $B_M(w)$ to the boundary of $\bar{B}_{M}(w)$.
By Lemma 2, there are $C_i=C_i(F)$ for $i=1,2$ such that
$${\bf P}[{\cal E}_w]\leq C_1M\exp(-C_2M).\eqno{(3.4)}$$
Note  that ${\cal E}_w$ and ${\cal E}_u$ are independent with the same distribution for fixed $w$ and $u$  if
$B_w(M)$ and $B_u(M)$ are two center squares of two different $3M$-squares.
With these observations and (2.1), 
if we take $M$ large, for small $C>0$, there exist $C_i=C_i(F)$ for $i=1,2$ in (3.4) and $C_j=C_j(F, M, C)$ for $j=3,4$
such that
\begin{eqnarray*}
&&{\bf P}[ \exists \mbox{ a NE path $\gamma$ from the origin with $|\gamma |\geq r$ and with less than $Cr$ closed edges}]\\
&\leq &\sum_{k\geq \lfloor r/ (M^2) \rfloor} 7^{2k} 2^k {\bf P}[ {\cal E}_{\bf 0}]^{k/30}\\
&\leq &\sum_{k\geq \lfloor r/ (M^2) \rfloor} 7^{2k} 2^k \left(C_1M \exp(-C_2 M)\right)^{k/30}\\
& \leq &C_3 \exp(-C_4 r).\hskip 4.6in  (3.5)
\end{eqnarray*}

{\bf Proof of Theorem 8.} On $\{\vec{T}({\bf 0}, (r, \theta)) \leq \epsilon^2(r-1)\}$,
 there exists a NE path $\gamma$ from the origin with 
$$|\gamma| \geq r-1 \mbox{ and }\vec{T}(\gamma) \leq \epsilon^2(r-1)$$
for some $\epsilon >0$. 
Note that if $|\gamma|\geq r-1$ and $\vec{T}(\gamma) \leq \epsilon^2 (r-1)$, then $\gamma$ contains  less than 
$\epsilon (r-1)$ closed edges.
So, if we take $C$ in (3.5) such that $C = \epsilon$ for some small $\epsilon$, by (3.5), there 
exist constants $C_i=C_i(F(0), \epsilon)$ for $i=2,3$  such that
\begin{eqnarray*}
&& {\bf P}[\vec{T}({\bf 0}, (r, \theta)) \leq \epsilon^2(r-1)]\\
&\leq &{\bf P}[ \exists \mbox{ a NE path $\gamma$ from the origin with $|\gamma|\geq (r-1)$, but $\vec{T}(\gamma) \leq \epsilon^2 (r-1)$ }]\\
&\leq &{\bf P}[ \exists \mbox{ a NE path $\gamma$ from the origin with $|\gamma |\geq (r-1)$ and with less than $C(r-1)$ closed edges}]\\
&\leq &C_2 \exp\left(-C_3 r\right).\hskip 4.7in {(3.6)}
\end{eqnarray*}
Therefore, for some $\delta >0$, there exist $C_i(F, \delta)$ for $i=1,2 $ such that
$${\bf P}[\vec{T}({\bf 0}, (r, \theta)) \leq \delta r]\leq C_1\exp(-C_2 r).\eqno{(3.7)}$$
Thus, Theorem 8 follows from (3.7). $\Box$

\section{Outside the percolation cone.}
The proofs for theorems outside the percolation cone also need the method of renormalization. 
We  assume that in Theorems 2 and 4
$$F(0)\geq  \vec{p}_c.\eqno{(4.1)}$$
Thus,  we say an edge is open or closed if $t(e)=0$ or $t(e) > 0$. With (4.1), we have
$${\bf P}[e \mbox{ is open}] \geq  \vec{p}_c.\eqno{(4.2)}$$
Similar to the proof in the last section, we work on a path $\gamma$ from the origin to $(r, \theta)$
and denote the $M$-squares of $\gamma$ by $\gamma_M$
for a large $M$.
If a square in $\gamma_M$ contains an edge $e$ with $t(e) > 0$, we say the square is a {\em bad} square.
Otherwise, it is a {\em good} square. 

If  there is a path $\gamma$ from the origin to $(r, \theta)$
with less than $Cr$ closed edges
for some small $C$, 
then there is less than $C|\gamma_M|$ bad squares in $\gamma_M$. 
Now we account the number of the choices for these squares. Note that $\gamma$ is connected, and so is $\gamma_M$.
We assume that
$$|\gamma_M|=k.\eqno{(4.3)}$$
As we proved in section 3, there are at most $7^{2k}$ choices for all
possible $\gamma_M$. When $\gamma_M$ is fixed, we select these bad squares. There are at most
$${k\choose Ck}\leq 2^k\eqno{(4.4)}$$ 
choices for these bad squares. 
We list all the bad squares as
$$S_1, S_2,\cdots , S_l$$
for $l\leq Ck$.
For each $S_i$,  path $\gamma$ will meet the boundary of $S_i$ at $v_i'$ and then use less than $2M$
edges to meet $v_i''$, another boundary point of $S_i$.
We denote the path from the origin to $v_1'$ by $\gamma_0$, from $v_1''$ to $v_2'$ by
$\gamma_1$, $\cdots$, from $v_{l-1}''$ to $v_{l}'$ by $\gamma_{l}$.
Note that bad edges are only contained in bad squares, so $\gamma_i$ does not contain a bad edge.
In other words, $\gamma_i$ is an open path for $i=1, \cdots, l$.
Now we reconstruct a NE fixed open path from $v_i'$ to $v_i''$. Let ${\cal E}_i$ be the event.
Since there exists a NE path with less than $2M$ edges from $v_i'$ to $v_i''$, 
$${\bf P}[{\cal E}_i]\geq (\vec{p}_c)^{2M}.\eqno{(4.5)}$$
Also, there are at most $(4M)^2$ choices for $v_i'$ and $v_i''$ when $S_i$ is fixed.
For a fixed $S_i$ and fixed $v_i'$ and $v_i''$ for $i=1,2,\cdots,l$,
$$\{\bigcap_{i=1}^l\{\exists \mbox{ open }\gamma_i\}\}\mbox{ and } \{\bigcap_{i=1}^l {\cal E}_i\}\mbox{ are independent,}\eqno{(4.6)}$$
since they use the outside and the inside edges of $S_i$, respectively.
With these reconstructions, on $\{\bigcap_{i=1}^l\{\exists \mbox{ open }\gamma_i\}\} \cap {\cal E}_i\}$, there exists
a NE open path from the origin to $(r, \theta)$.

If we assume that $|\gamma_M|=k,$ by Lemma 1, there are at least $k/15$ disjoint center $3M$-squares of $\gamma$.
For each center $3M$-square, it contains at least $M$ vertices of $\gamma$. Thus, 
$$|\gamma|\geq M |\gamma_M|/15= kM /15.\eqno{}$$ 
Note also that $\gamma$ is directed, so together the  inequality above,
$$2r \geq |\gamma|\geq  kM /15.\eqno{(4.7)}$$

With these observations, we have
\begin{eqnarray*}
&&{\bf P}[ \exists \mbox{ a NE $\gamma$ from ${\bf 0}$ to $(r, \theta)$ with less than $Cr$ closed edges}]\\
&\leq &\sum_{k\geq \lfloor r/ M^2\rfloor} 7^{2k} 2^k \sum_{l=1}^{Ck}(4M)^{2 l}{\bf P}\left[\bigcap_{i=1}^l\{\exists \mbox{ open }\gamma_i\}\right]\\
&\leq &\sum_{k\geq \lfloor r/ M^2\rfloor} 7^{2k} 2^k \sum_{l=1}^{Ck}(4M)^{2 l}{\bf P}\left[\bigcap_{i=1}^l\{\exists \mbox{ open }\gamma_i\}\right](\vec{p}_c)^{-2CkM}\prod_{i=1}^l {\bf P}\left[{\cal E}_i\right]\\
&\leq & \sum_{k\geq \lfloor r/ M^2\rfloor} k7^{2k} 2^k(4M)^{2 Ck}(\vec{p}_c)^{-2CkM}{\bf P}\left[\bigcap_{i=1}^l\{\exists \gamma_i\}\cap {\cal E}_i\right]\\
&\leq & \sum_{k\geq \lfloor r/ M^2\rfloor} k7^{2k} 2^k (4M)^{2 Ck}(\vec{p}_c)^{-2CkM}{\bf P}[\exists \mbox{ an open path
from the origin to $(r, \theta)$}].\hskip 1cm (4.8)
\end{eqnarray*}
 Let $u=(u_1, u_2)$  be the ending vertex of $\gamma$. If $\theta< \theta_p^-$, then
$$sl ({\bf 0}, u)=\tan(\theta) < \tan(\theta_p^-).\eqno{(4.9)}$$
In addition, 
$$u_1= O(r).\eqno{(4.10)}$$
Thus by Lemma 3,  (4.9), (4.10), and (4.7), there exist $C_i=C_i(F, \theta)$ for $i=1,2,3$ such that
$${\bf P}[\exists \mbox{ an open path from the origin to $(r, \theta)$}]\leq C_1 \exp(-C_2 r)\leq C_1\exp(-C_3 Mk). \eqno{(4.11)}$$

If we substitute (4.11) into (4.8), there exist $C_i=C_i(F, \theta)$ for $i=1,2$ such that
\begin{eqnarray*}
&&{\bf P}[ \exists \mbox{ a NE $\gamma$ from ${\bf 0}$ to $(r, \theta)$ with less than $Cr$ closed edges}]\\
&\leq &\sum_{k\geq \lfloor r/ M^2\rfloor} k7^{2k} 2^k k(4M)^{2 Ck}(\vec{p}_c)^{-2CkM} C_1\exp(-C_2 Mk/15).
\end{eqnarray*}
Therefore, if we take $M$ large and then $C=C(M, C_2)$ small, there exist $C_i=C_i(F,\theta, C)$ for
$i=3,4$ such that
$${\bf P}[ \exists \mbox{ a NE $\gamma$ from ${\bf 0}$ to $(r, \theta)$ with less than $Cr$ closed edges}]\leq C_3
\exp(-C_4 r).$$
In summary, if $\theta< \theta_p^-$, for all $r$, there exist $C=C(F, \theta)$ and $C_i=C_i(F, \theta, C)$
for $i=1,2$ such that
$${\bf P}[ \exists \mbox{ a NE $\gamma$ from ${\bf 0}$ to $(r, \theta)$ with less than $Cr$ closed edges}]\leq C_1
\exp(-C_2 r).\eqno{(4.12)}$$
With (4.12), we show Theorems 2 and 4:\\

{\bf Proofs of Theorems 2 and  4.}
Suppose that there exists a NE path $\gamma$ from the origin to $(r, \theta)$ with
$$\vec{T}(\gamma) \leq C_1 r.$$
By (4.12), we may choose  the $C$ in (4.12) such that 
\begin{eqnarray*}
&&{\bf P}[\vec{T}(\gamma) \leq C_1 r]\\
&=&{\bf P}[ \vec{T}(\gamma) \leq C_1 r,  \mbox{$\gamma$  with  more than $Cr$ closed edges}]\\
&&+{\bf P}[ \vec{T}(\gamma) \leq C_1 r, \mbox{$\gamma$  with  less than $Cr$ closed edges}]\\
&\leq & {\bf P}[ \vec{T}(\gamma) \leq C_1 r, \mbox{$\gamma$  with  more than $Cr$ closed edges}]+C_2\exp(-C_3 r).\hskip 3cm (4.13)
\end{eqnarray*}
For each closed edge $e$, we know that $t(e) > 0$. For $\epsilon >0$, we  take $\delta >0$ small such that
$${\bf P}[0< t(e) \leq \delta ] =F(\delta)-F(0)\leq \epsilon.\eqno{}$$
For each closed edge, if it satisfies $t(e) \leq \delta$, we say it is a bad edge. Thus,
$${\bf P}[ e \mbox{ is closed and bad}] ={\bf P}[0< t(e) \leq \delta]\leq \epsilon.\eqno{}$$
Now, on  $\{\exists$ a NE $\gamma$ from ${\bf 0}$ to $(r, \theta)$ with more than $C r$ closed edges$\}$,
we estimate the event that there are at least $Cr/2$ bad edges in $\gamma$. By (4.7),
$$|\gamma|\leq 2r.\eqno{}$$
Now we fix path $\gamma$. Since each vertex in $\gamma$ can be adjacent only from a north or an east edge,
there are at most $2^{2r}$ choices for  $\gamma$.
If $\gamma$ is fixed, there are at most 
$$\sum_{l=1}^{2r} {2r \choose l}\leq  2^{2r}$$
choices for these closed edges. If these closed edges are fixed, as we mentioned above, each edge has a probability
less than $\epsilon$ such that it  is also bad. In addition, we also have another $2^{2r}$ choices to select
these bad  edges from these closed edges.
Therefore,  if we take $\epsilon=\epsilon(F, \delta, C)$ small,  then there exist $C_i=C_i(F, \delta, C)$ for $i=2,3$
such that
\begin{eqnarray*}
&&{\bf P}[ \exists \mbox{ a NE $\gamma$ from ${\bf 0}$ to $(r, \theta)$ with more than $Cr$ closed edges}, \\
&&\hskip 1cm \mbox{ these closed edges contain more than $Cr/2$   bad edges}]\\
&\leq &\sum_{l=Cr/2}^\infty 2^{2r} 2^{2r} 2^{2r} (\epsilon)^l\\
&\leq &C_2 \exp(-C_3 r). \hskip 4.5in (4.14)
\end{eqnarray*}
If
$$\vec{T}({\bf 0}, (r, \theta))\leq C_1r$$
for $\theta < \theta_p^-$,
then there is a NE path $\gamma$ from ${\bf 0}$ to $(r, \theta)$ with a passage time less than $C_1r$.
Therefore, by (4.13) and (4.14),
\begin{eqnarray*}
&&{\bf P}[ \vec{T}({\bf 0}, (r, \theta))\leq C_1r]\\
&\leq &{\bf P}[ \exists \mbox{ a NE $\gamma$ from ${\bf 0}$ to $(r, \theta)$ with more than $C r$ closed edges}, \\
&&\hskip 0.3cm \mbox{ these closed edges contain less than $Cr/2$   bad edges}, \vec{T}(\gamma) \leq C_1r]+C_2 \exp(-C_3 r)\\
&\leq &{\bf P}[ \exists \mbox{ a NE $\gamma$ from ${\bf 0}$ to $(r,\theta)$, $\gamma$ contains less than $Cr/2$   bad edges}, \vec{T}(\gamma) \leq C_1r]\\
&&+C_2 \exp(-C_3 r).
\hskip 4.5in (4.15)
\end{eqnarray*}
If there is a NE path from ${\bf 0}$ to $(r, \theta)$ with  less than $Cr/2$   bad edges among these $Cr$ closed edges, 
note that each good edge costs at least passage time $\delta$, so
 the passage time of the path is more than 
$\delta Cr/2.$
Thus, if we select $C_1$ such that
$$ C_1 <C \delta /2,$$
$$\!\!\!\!\!\!{\bf P}[ \exists \mbox{ a NE $\gamma$ from ${\bf 0}$ to $(r,\theta)$, $\gamma$ contains less than $Cr/2$   bad edges}, \vec{T}(\gamma) \leq C_1r]=0.\eqno{(4.16)}$$
By (4.15) and (4.16), for $F(0) \geq \vec{p}_c$, $\theta< \theta_p^-$,  there exist $C_1=C_1(F, \theta)$
and $C_i=C_i(F, \theta, C_1)$ for $i=2,3$ such that
$${\bf P}[\vec{T}({\bf 0}, (r, \theta)) \leq C_1 r]\leq C_2 \exp(-C_3 r).\eqno{(4.17)}$$
When $\theta > \theta_p^+$, by symmetry,  we still have (4.17).
Therefore, Theorems 2 and 4 follow. $\Box$

\section{ Inside the percolation cone.}
In section 5, we assume that $F(0) =p > \vec{p}_c$ and $\theta \in [\theta_p^-,   \theta_p^+]$.
Edge $e$ is called an open or a closed edge if $t(e) =0$ or $t(e) >0$, respectively.
We define $\tau(e)=0$ if $t(e)=0$, or $\tau(e)=1$ if $t(e) >0$.
We also denote by the passage time $\vec{T}_\tau(u, v)$ corresponding to $\tau(e)$.
Let
$$B_\tau(t)= \{v\in {\bf Z}^2: \vec{T}_\tau({\bf 0}, v)\leq t\}.\eqno{(5.1)}$$
We also assume that $(r, \theta)\in {\bf Z}^2$ without loss generality.
If $(r, \theta)\in B_\tau(t)$, then
$$\vec{T}_\tau({\bf 0}, (r, \theta))\leq t.\eqno{(5.2)}$$
Note that $B_\tau(t)$ will eventually cover all the vertices in  ${\bf R}\times [0, \pi/2]$ as $t\rightarrow \infty$, so there exists
a $t$ such that $(r, \theta)\in B_\tau(t)$. Let
$\sigma$ be the smallest $t$ such that $(r, \theta)\in B_\tau(t)$.
We will estimate $\sigma$ to show that there exist $C_i=C_i(F)$ for $i=1,2$ such that for all large $k$,
$${\bf P}[\sigma \geq k] \leq C_1 \exp(-C_2 k).\eqno{(5.3)}$$
Note that
$${\bf P}[\sigma \geq k]=\sum_{\Gamma }{\bf P}[\sigma\geq k, B_\tau (k-2) =\Gamma],$$
where $\Gamma$, containing the origin,   takes all possible vertex sets in the first quadrant.
We also remark that for $\Gamma_1$ and $\Gamma_2$,
$$\{\sigma\geq k, B_\tau (k-2) =\Gamma_1\}\mbox{ and } \{\sigma\geq k, B_\tau (k-2) =\Gamma_2\}\mbox{ are disjoint}.$$
If $\sigma\geq k$ and $B_\tau (k-2) =\Gamma$, then $\Gamma$ does not contain  $(r, \theta)$:
$$\Gamma\cap (r,\theta)=\emptyset.\eqno{(5.4)}$$
In other words,  all $\Gamma$ in the above sum do not contain $(r, \theta)$.
Thus, by Lemma 5, there is no  NE open path from $\partial_o(\Gamma)$ to $(r, \theta)$, without using edges of $\Gamma\cup \partial_e \Gamma$.
Otherwise,
$\sigma < k$, which is contrary to the assumption that $\sigma \geq k$.
For a fixed $\Gamma$, we denote by ${\cal E}_{k}(\Gamma)$ the above event that 
there is no  NE open path without using edges of $\Gamma\cup \partial_e \Gamma$
from $\partial_o(\Gamma)$ to $(r, \theta)$.
Note that ${\cal E}_k(\Gamma)$ only depends on configurations of edges outside $\Gamma\cup \partial_e \Gamma$, 
so by Lemma 6, 
for any fixed $\Gamma$ with $\Gamma\cap (r,\theta)=\emptyset$,
$$ {\cal E}_{k}(\Gamma) \mbox{ and } \{B_\tau (k-2) =\Gamma \}\mbox{ are independent}.\eqno{(5.5)}$$
Note that if there is a NE open path from the origin to $(r, \theta)$, then
there exists a NE open path outside $\Gamma\cup \partial_e \Gamma$ from $\partial_o\Gamma$ to $(r, \theta)$. 
By (1.5), there exists $0<\delta<1$ such that
for a fixed $\Gamma$,
$${\bf P}[{\cal E}_{k}(\Gamma)]\leq 1-{\bf P}[{\bf 0} \rightarrow (r, \theta)]\leq 1-\delta .\eqno{(5.6)}$$

With these observations, 
\begin{eqnarray*}
{\bf P}[\sigma \geq k]
&=&\sum_{\Gamma }{\bf P}[\sigma\geq k, B_\tau (k-2) =\Gamma ]\\
&\leq &\sum_{\Gamma }{\bf P}[B_\tau (k-2) =\Gamma, (r, \theta)\not\in \Gamma,  {\cal E}_{k}(\Gamma)]\\
&\leq&\sum_{\Gamma }{\bf P}[B_\tau (k-2) =\Gamma, (r, \theta)\not\in \Gamma](1-\delta).
\end{eqnarray*}
Note that for a fixed $\Gamma$, by Lemma 5 again,
$$\{B_\tau (k-2) =\Gamma, (r, \theta)\not\in \Gamma\}\subset \{\sigma \geq (k-2)\}.$$
Therefore, 
$${\bf P}[\sigma \geq k] \leq (1-\delta) {\bf P}[\sigma \geq (k-2)].\eqno{(5.7)}$$
Thus, (5.3) follows if we iterate (5.7). We show Theorem 1 by (5.3). In fact, if $t(e)$ is bounded from above by a constant,
then Theorem 1 is implied by (5.3) directly. However, if we restrict ourselves only on a moment condition,
the proof is complicated, as follows:\\

{\bf Proof of Theorem 1.}
On $\{\sigma = k\}$, there exists an optimal  path $\gamma_k$ in $\tau(e)$  from the origin to $(r, \theta)$
with only $k$ edges $\{e_i\}$ such that $\tau(e_i)=1$. 
For each configuration on $\{\sigma =k\}$, we use a unique way to select an optimal path in $\tau(e)$  with
these $k$  edges. We still denote the path by $\gamma_k$. For configuration $\omega$, and    path $\gamma_k(\omega)$,
let $e_1, e_2, \cdots, e_k \subset \gamma_k$ with $\tau(e_i) =1$.
 Note that on $\{\sigma = k\}$, only $t(e_i) >0$, but the others are zero-edges, so
\begin{eqnarray*}
&&{\bf E}\left(\vec{T}({\bf 0}, (r, \theta)\right)^m=\sum_{k=1}^\infty {\bf E}\left[\left(\vec{T}({\bf 0}, (r, \theta)\right)^m ; \sigma = k\right]\\
&=&\sum_{k=1}^\infty\sum_{\beta}\sum_{e_i, 1\leq i\leq k} {\bf E}\left[\left(\sum_{i=1}^k t(e_i)\right)^m ; \gamma_k =\beta, 
\tau(e_i) =1, \tau(e_s)=0, s\neq i, 1\leq i\leq k\right],\\
\end{eqnarray*}
where $\beta$ is a fixed NE path, the second sum takes over all possible NE paths $\beta$ that are
from the origin to $(r, \theta)$, and the third sum takes over all possible $k$ edges $e_i$ for $i=1,2,\cdots, k$ on path $\beta$.
With these decompositions,  the events
$$\{\gamma_k =\beta, \tau(e_i) =1, \tau(e_s)=0, s\neq i,1\leq i\leq k \}\mbox{ are disjoint} \eqno{(5.8)}$$
for different paths $\beta$ and different selections of $e_i$ in $\beta$.

On the event in (5.8) for a fixed $\beta$ and  these fixed $e_i$, $1\leq i\leq k$, we denote by  event ${\cal E}_j$ if 
$$t(e_j)=\max \{t(e_1), \cdots , t(e_k)\} \mbox{ and } t(e_i)< t(e_j) \mbox{ for } i=1,\cdots, j-1.$$
Note that ${\cal E}_j$ are disjoint, so
\begin{eqnarray*}
&&{\bf E}\left[\vec{T}({\bf 0}, (r, \theta))\right]^m\\
&\leq & \sum_{k=1}^\infty\sum_{\beta}\sum_{e_i, 1\leq i\leq k}\sum_{j=1}^k{\bf E}\left[\left( k t(e_j)\right)^m ; \gamma_k =\beta,\tau(e_i) =1, \tau(e_s)=0, s\neq i,1\leq i\leq k, {\cal E}_j\right]\\
&\leq & \sum_{k=1}^\infty\sum_{\beta}\sum_{e_i, 1\leq i\leq k}\sum_{j=1}^k{\bf E}[( k t(e_j))^m ; \gamma_k =\beta,
\tau(e_i) =1, \tau(e_s)=0, s\neq i,1\leq i\leq k].\hskip 1.5cm (5.9)
\end{eqnarray*}
Note that $\vec{T}({\bf 0}, (r, \theta))$ only depends on finite many edges. Note also that any distribution can be 
approximate by a discrete distribution (see Theorem 4.13 in Wheeden and  Zygmund (1977)). We may assume that
$t(e)$ is a discrete random variable taking values in an accountable set $\{l\}$ without loss of generality.
 Thus, for fixed $k$, $\beta$, and $j$,
\begin{eqnarray*}
&&{\bf E}\left[( k t(e_j))^m ; \gamma_k =\beta, \tau(e_i) =1, \tau(e_s)=0, s\neq i, 1\leq i\leq k \right]\\
=&&\sum_{l>0} (kl)^m {\bf P}[\gamma_k =\beta, \tau(e_i) =1, \tau(e_s)=0, s\neq i,1\leq i\leq k, i\neq j,   t(e_j)=l].\hskip 1.8cm (5.10)
\end{eqnarray*}
We will show that
\begin{eqnarray*}
&&{\bf P}[  \gamma_k =\beta, \tau(e_i) =1, \tau(e_s)=0, s\neq i,1\leq i\leq k, i\neq j,   t(e_j)=l ]\\
&\leq & (1-F(0))^{-1}{\bf P}[\gamma_k =\beta, \tau(e_i) =1, \tau(e_s)=0, s\neq i,1\leq i\leq k]  {\bf P}[ t(e_j) =l].
\hskip 2cm (5.11)
\end{eqnarray*}
For each configuration $\omega$ with 
$$\omega \in\{ \gamma_k =\beta, \tau(e_i) =1, \tau(e_s)=0, s\neq i,1\leq i\leq k, i\neq j,  t(e_j)=l\},$$
we assume that  
$$\omega=(\omega(e_j)=l, \omega'(e_j)),\eqno{}$$
where
$\omega'(e_j)$ is the restriction of $\omega$ on ${\bf Z}^2\setminus e_j$.
Let 
$$W'_j=\{\omega'(e_j): (\omega(e_j)=l, \omega'(e_j) )\in \{ \gamma_k =\beta, \tau(e_i) =1, \tau(e_s)=0, s\neq i,1\leq i\leq k, i\neq j,  t(e_j)=l\}\}.$$
Since $(t(e_j), \tau(e_j))$ and $(t(e_i), \tau(e_i))$ are independent when $i\neq j$,
\begin{eqnarray*}
&&{\bf P}[  \gamma_k =\beta, \tau(e_i) =1, \tau(e_s)=0, s\neq i,1\leq i\leq k, i\neq j,   t(e_j)=l ]\\
&=&\sum_{\omega'(e_j)\in W'(e_j)}{\bf P}[  \{t(e_i)\}_{i\neq j}=\omega'(e_j),t(e_j)=l ]\\
&=& \sum_{\omega'(e_j)\in W'(e_i)}{\bf P}[ \{t(e_i)\}_{i\neq j}=\omega'(e_j)]{\bf P}[t(e_j)=l]\left({{\bf P}[t(e_j)>0]\over {\bf P}[t(e_j)>0]}\right)\\
&= &\sum_{\omega'(e_j)\in W'(e_i)}{\bf P}[  \{t(e_i)\}_{i\neq j}=\omega'(e_j), t(e_j) > 0]{\bf P}[t(e_j)=l] (1-F(0))^{-1}\\
&= & {\bf P}[t(e_j)=l] (1-F(0))^{-1}
\sum_{\omega'(e_j)\in W'(e_i)}{\bf P}[  \{t(e_i)\}_{i\neq j}=\omega'(e_j), t(e_j) > 0].
\end{eqnarray*}
Note that for different $\omega'_1(e_j)$ and $\omega'_2(e_j)$,
$$\{ \{t(e_i)\}_{i\neq j}=\omega'_1(e_j), t(e_j)>0\} \mbox{ and } \{ \{t(e_i)\}_{i\neq j}=\omega'_2(e_j), t(e_j)>0\}\mbox{ are disjoint}.$$
Note also that  $t(e_j)>0$ implies that $\tau(e_j)=1$ and the configurations of the other edges keep the same  in $W_j'$, so
$$\{\{t(e_i)\}_{i\neq j}=\omega'(e_j), t(e_j) > 0\}\subset \{\gamma_k =\beta, \tau(e_i) =1, \tau(e_s)=0, s\neq i,1\leq i\leq k\}.$$
With these observations,
\begin{eqnarray*}
&&{\bf P}[  \gamma_k =\beta, \tau(e_i) =1, \tau(e_s)=0, s\neq i,1\leq i\leq k, i\neq j,   t(e_j)=l ]\\
&\leq &(1-F(0))^{-1}{\bf P}[t(e_j) =l]{\bf P}[\gamma_k =\beta, \tau(e_i) =1, \tau(e_s)=0, s\neq i, 1\leq i\leq k].
\end{eqnarray*}

Therefore, (5.11) follows.
If we apply (5.11) in (5.10), by the assumption that $Et^m(e) < \infty$, there exist $C=C(F)$ and $C_1=C_1(F,m)$
such that
\begin{eqnarray*}
&&{\bf E}[( k t(e_j))^m ; \gamma_k =\beta, t(e_i) >0, t(e_s)=0, s\neq i,1\leq i\leq k]\\
&\leq & C k^m \sum_{l > 0} l^m {\bf P}[t(e_j)=l]{\bf P}[\gamma_k =\beta, \tau(e_i) =1, \tau(e_s)=0, s\neq i,1\leq i\leq k]\\
&\leq &C_1 k^m {\bf P}[\gamma_k =\beta, \tau(e_i) =1, \tau(e_s)=0, s\neq i,1\leq i\leq k ]. \hskip 4.5cm {(5.12)}
\end{eqnarray*} 
We apply (5.12) in (5.9):
\begin{eqnarray*}
&&{\bf E}\left(\vec{T}({\bf 0}, (r, \theta))\right)^m\\
&\leq &C \sum_{k=1}^\infty\sum_{\beta}\sum_{e_i, 1\leq i\leq k}\sum_{j=1}^k k^{m} {\bf P}[\gamma_k =\beta, \tau(e_i) =1, \tau(e_s)=0, s\neq i,1\leq i\leq k ]\\
&\leq &C \sum_{k=1}^\infty k^{m+1}\sum_{\beta}\sum_{e_i, 1\leq i\leq k} {\bf P}[\gamma_k =\beta, \tau(e_i) =1, \tau(e_s)=0, s\neq i,1\leq i\leq k ]. \hskip 2cm {(5.13)}
\end{eqnarray*}
Note that event
$$\{\gamma_k =\beta, \tau(e_i) =1, \tau(e_s)=0, s\neq i,1\leq i\leq k \}$$
implies that an optimal path in $\tau(e)$ from the origin to $(r, \theta)$
contains  only $k$  1-edges.
Therefore, it implies that $\{\sigma \geq k\}$.
With these observations, by  the disjoint property of (5.8) and (5.13),
$${\bf E}\left(\vec{T}({\bf 0}, (r, \theta))\right)^m\leq C\sum_{k=1}^\infty k^{m+1}{\bf P}[ \sigma \geq k]. \eqno{ (5.14)}$$

Finally, by (5.3) and (5.14), for all $r$ and all $\theta\in [\theta_p^-, \theta_p^+]$,
there exists $C_3=C_3(F, m)$ such that 
$${\bf E}\left(\vec{T}({\bf 0}, (r, \theta))\right)^m\leq C \sum_{k=1}^\infty k^{m+1}{\bf P}[ \sigma \geq k]\leq C \sum_{k=1}^\infty k^{m+1}C_1\exp(-C_2 k) \leq C_3. \eqno{(5.15)}$$
Theorem 1 follows from (5.15). $\Box$

\section{ Critical phase.} 
 
{\bf Proof of Theorem 5.} First we show (1.7) in Theorem 5. We may take $h$ small 
such that $F$ has two different situations at $F(0)$:

(a) $F(x)=F(0)=\vec{p}_c $ is  a constant on $[0, h]$.

(b) $\exists $ a sequence $\{x_n\}$ with $x_n \downarrow 0$ such that $F(x_n)\downarrow F(0)$.\\
Let us assume that case (b) holds. For all $\epsilon \leq h$, we construct another distribution:
$$G_n (x)=\left\{\begin{array}{cc}
F(0) &\mbox{ if  $0\leq x< x_n$,}\\
F(x)&\mbox{ if  $x_n\leq x.$}
\end{array}
\right.
$$
By this definition, for each $n$,
$$G_n (x_n) > \vec{p}_c.$$
By (1.5), for all $r$,  there exists a directed path from the origin to 
$(r, \pi/4)$ such that
its passage time in each edge is at most $x_n$ with a positive probability. By (1.1), for each $n$,
$$\vec{\mu}_{G_n}(\pi/4)\leq 2x_n.\eqno{(6.1)}$$
By Lemma 7, 
$$\vec{\mu}_F(\pi/4)\leq \vec{\mu}_{G_n}(\pi/4)\leq 2x_n. \eqno{(6.2)}$$
By (6.2), we can show that
$$\vec{\mu}_{F}(\pi/4)=0.\eqno{(6.3)}$$
Therefore, (1.7) in Theorem 5 follows if case (b) holds.

Now we focus on case (a). 
Note that $F$ cannot be flat forever, so there are points $h_1 > h >0$ such that $F(h_1) > F(0)$ and $F(x)=F(0)$
for $0\leq x \leq h$.
Now we assume that $F$ satisfies the following extra condition:

(i). There exists $h>0$ such that
$$F(x)= F(0)=\vec{p}_c, \mbox{ when } 0\leq x< h, \mbox{ and } F(x) > \vec{p}_c,\mbox{ when } x\geq h.$$
In other words, there is a jump point at $h$.

We focus on case (a) (i). We take $\epsilon>0$ small such that
$$F(0) +\epsilon < F(h).$$
Then we construct another distribution:
$$G_\epsilon (x)=\left\{\begin{array}{cc}
F(0)+\epsilon &\mbox{ if  $0\leq x < h$,}\\
F(h)&\mbox{ if $x=h$,}\\
F(x)&\mbox{ if  $h< x.$}
\end{array}
\right.
$$
As we defined,  $t(e)$ is the random variable with distribution $F$. Let
$g_\epsilon(e)$ be the random variable with distribution $G_\epsilon$. In addition, 
$t(e)$ and $g_\epsilon(e)$ are signed values independently edge by edge as we defined.
The key step is to couple these two random variables together.
Define $g_\epsilon(e)$ as follows:\\
If $t(e)=0$, then $g_\epsilon(e)=0$.\\
If $t(e)=x> h$, then  $g_\epsilon(e)=x$.\\
If $t(e)=h$, then 
$$
g_\epsilon (e)\left\{\begin{array}{cc}
 =0&\mbox{ with probability $\epsilon \left(F(h)-\vec{p}_c\right)^{-1}$,}\\
 >0 &\mbox{ with  probability $1-\epsilon \left(F(h)-\vec{p}_c\right)^{-1}.$}
\end{array}
\right.
$$
Now we verify that $g_\epsilon$ has  distribution $G_\epsilon$. 
For short, we simply replace  $g_\epsilon(e)$ and $t(e)$ by $g_\epsilon$ and $t$:
\begin{eqnarray*}
&&{\bf P}[ g_\epsilon =0]\\
&=&{\bf P}[ g_\epsilon =0\,\, |\,\, t=0] \vec{p}_c +{\bf P}[ g_\epsilon =0\,\, |\,\, t=h] (F(h)-\vec{p}_c)+{\bf P}[ g_\epsilon =0\,\, |\,\, t>h]{\bf P}[ t > h]\\
&=& \vec{p}_c+\epsilon =F(0)+\epsilon .\hskip 4.5in (6.4)
\end{eqnarray*}
Note that 
$$ {\bf P}[ g_\epsilon > h]={\bf P}[ t > h],$$
so by (6.4),
$${\bf P}[ g_\epsilon =h]=1-{\bf P}[ g_\epsilon =0]-{\bf P}[ g_\epsilon > h]=1-(\vec{p}_c+\epsilon)-(1-F(h))=F(h)-(\vec{p}_c+\epsilon).\eqno{(6.5)}$$
Finally, for $x > h$, by (6.4) and (6.5), 
$${\bf P}[ g_\epsilon \leq  x]={\bf P}[ h< g_\epsilon \leq  x]-{\bf P}[ g_\epsilon \leq  h]=
{\bf P}[ h< t \leq  x]-{\bf P}[ g_\epsilon \leq  h]=F(x).\eqno{(6.6)}$$
Thus, by (6.4)--(6.6), $g_\epsilon$ indeed has  distribution $G_\epsilon$. 

Now we show  (1.7) in Theorem 5 under case (a) (i). 
Let $\gamma_t$ be an optimal path for $\vec{T}_{t}({\bf 0}, (r, \pi/4))$ with time state $t(e)$, and let
$\gamma_{g_\epsilon}$ be an optimal path for $T_{g_\epsilon}({\bf 0}, (r, \pi/4))$  with time state $g_\epsilon(e)$.
Here, for each configuration, we select $\gamma^{g_\epsilon}$ in a unique method.
For each edge $e\in \gamma_{g_\epsilon}$, we consider passage time $t(e)$.
If $t(e)> h$, then $g_\epsilon(e)=t(e)$, as we defined. 
If $t(e)=0$, then   $g_\epsilon(e)=0$. In addition, if $t(e)=h$, it follows from the definition that  $g_\epsilon(e)\leq h$.
Therefore, 
$$\vec{T}_t({\bf 0}, (r, \pi/4))\leq \vec{T}(\gamma_{g_\epsilon}) + h\sum_{e\in \gamma_{g_\epsilon}} I_{(t(e)=0, g_\epsilon(e)=h)}.\eqno{(6.7)}$$
By (6.7),
$${\bf E}\vec{T}_t({\bf 0},(r, \pi/4))\leq {\bf E}\vec{T}(\gamma_{g_\epsilon}) + h\sum_{\beta}\sum_{e\in \beta } {\bf P}[t(e)=0, g_\epsilon(e) =h ,\gamma_{g\epsilon}=\beta],\eqno{(6.8)}$$
where the first sum in (6.8) takes over all possible NE paths $\beta$ from ${\bf 0}$ to 
 $(r, \pi/4)$.
Let us estimate 
$$\sum_{\beta}\sum_{e\in \beta } {\bf P}[t(e)=0, g_\epsilon(e) =h ,\gamma_{g_\epsilon}=\beta].$$
Note that the value of $g_\epsilon(e)$ may depend on the value of $t(e)$, but not on the other values of $t(b)$ for $b\neq e$, so
by our definition,
\begin{eqnarray*}
&&{\bf P}[t(e)=0, g_\epsilon(e)=h,\gamma_{g_\epsilon}=\beta]={\bf P}[g_\epsilon(e)=h\,\, |\,\, t(e)=0, \gamma_{g_\epsilon}=\beta]{\bf P}[t(e)=0,\gamma_{g_\epsilon}=\beta]\\
&\leq & {\bf P}[g_\epsilon(e)=h\,\, |\,\, t(e)=0]{\bf P}[\gamma_{\epsilon}=\beta]=\epsilon (F(h)-\vec{p}_c)^{-1}{\bf P}[\gamma_{g_\epsilon}=\beta].\hskip 4cm (6.9)
\end{eqnarray*}
By   (6.9), note that $\beta$ has at most $Cr$ edges, so
$$\!\!\!\sum_{\beta}\sum_{e\in \beta } {\bf P}[t(e)=0, g_\epsilon(e)=h,\gamma^t=\beta]\leq \sum_{\beta}\sum_{e\in \beta } \epsilon(1-\vec{p}_c)^{-1}{\bf P}[\gamma^t=\beta]\leq 2 \epsilon (F(h)-\vec{p}_c)^{-1}r.\eqno{(6.10)}$$
By  (6.8) and (6.10), there exists $C=C(F)$ such that
$${\bf E}{{\bf T}_t({\bf 0}, (r, \pi/4))\over r}\leq {\bf E} {\vec{T}_{g_\epsilon}({\bf 0}, (r, \pi/4))\over r}+C\epsilon.\eqno{(6.11)}$$
We take $r\rightarrow \infty$ in (6.11) to have
$$\vec{\mu}_F(\pi/4)\leq \vec{\mu}_{G_\epsilon}(\pi/4)+C\epsilon.\eqno{(6.12)}$$
Note that $G_\epsilon(0)> \vec{p}_c$, so by Corollary 3 and (6.12),
$$\vec{\mu}_F(\pi/4)=0.$$
Therefore, (1.7) in Theorem 5 follows under case  (a) (i).

Finally, we focus on case (a) without other assumptions.  As we mentioned, $t(e)$ is not a constant. Thus,
 there exists $h_1 > h$ such that $F(h_1) > F(0)$.
We construct
$$H(x)=\left\{\begin{array}{cc}
F(0) &\mbox{ if  $0\leq x < h_1$,}\\
F(x)&\mbox{ if $h_1\leq x$.}
\end{array}
\right.
$$
With this definition,
$$H \leq F, \mbox{ and }H(0)=\vec{p}_c,$$
and $H(x)$ has a jump point at $h_1$. By the analysis of case (a) (i), we have
$$\vec{\mu}_{H}(\pi)=0.\eqno{(6.13)}$$
By Lemma 7,
$$\vec{\mu}_{F}(\pi/4)\leq \vec{\mu}_{H}(\pi/4)=0.\eqno{(6.14)}$$
Thus, (1.7) in Theorem 5 under  case (a)  follows from (6.14).
If we put cases (a) and  (b) together,  (1.7) in Theorem 5 follows.\\

Now we show (1.8) in Theorem 5. When $F(0) =p > \vec{p}_c$, by Corollary 3, $\vec{\mu}_F(\pi/4)=0$.
Thus
$$\lim_{p\downarrow \vec{p}_c} \vec{\mu}_F(\pi/4)=0.\eqno{(6.15)}$$
If $F(0)=p < \vec{p}_c$, we take $\epsilon $ small such that
$$\vec{p}_c < p+\epsilon.\eqno{(6.16)}$$
Thus, we use the same method in case (a) to construct $G_\epsilon$ with $G_\epsilon(0) >\vec{p}_c$
and $G_\epsilon(x) \geq F(x)$.
By the same proof of (a) (i), we can show that $\mu_F(\pi/4)$ is bounded by $\mu_{G\epsilon}+C\epsilon$
from above. Thus, by using Corollary 3, we have 
$$0\leq \lim_{p \uparrow \vec{p}_c} \vec{\mu}_F(\pi/4)\leq \lim_{\epsilon \rightarrow 0} \mu_{G_\epsilon}(\pi/4)=0.\eqno{(6.17)}$$
Therefore, (1.8) in Theorem 5 follows. $\Box$\\

{\bf Proof of Theorem 7.} In this proof, we assume that $t(e)$ only takes 0 (open) and 1 (closed)
with probability $\vec{p}_c$
and $1-\vec{p}_c$, respectively. Let $L_r$ be the line $y=-x+r$ inside the first quadrant. Note that $L_0$ is just
the origin.
Bezuidenhout and Grimmett (1991) showed that for fixed $L_{r_1}$ and for $0< \delta <1$, there exists $r_2=r_2(r_1)$ such that
$${\bf P}[ L_{r_1} \not\rightarrow L_{r_2}] \geq \delta.\eqno{(6.18)}$$
If $L_{r_1} \not\rightarrow L_{r_2}$, then  NE path from $L_{r_1}$ to $L_{r_2}$
has to use at least one edge with passage time 1. Let $I(L_{r_1},L_{r_2})$ be the indicator of the event that
there is no  NE open path from $L_{r_1}$ to $L_{r_2}$.
For large $r$, let
$r_1\leq r_2\leq \cdots \leq r_m\leq r$ be a sequence  such that for $i< m-1$,
$${\bf P}[ L_{r_i} \not\rightarrow L_{r_{i+1}}] \geq \delta.\eqno{(6.19)}$$
Note that any NE  path from the origin to $(r, \pi/4)$ has to cross the strip between $L_i$ and $L_{i+1}$ 
for $i=1,\cdots, m-1$, so
$${\bf E}\vec{T}({\bf 0}, (r, \pi/4)) \geq {\bf E}\sum_{i=1}^m I(L_{r_1},L_{r_2})= \delta m.\eqno{(6.20)}$$
By (6.18), we have $m\rightarrow \infty$ as $r\rightarrow \infty$. Therefore, by (6.20),
$$\lim_{r\rightarrow \infty} {\bf E}\vec{T}({\bf 0}, (r, \pi/4))=\infty.\eqno{(6.21)}$$
Theorem 7 follows from (6.21). $\Box$

\newpage

\begin{center}
{\bf \large References}
\end{center}
Alexander, K. (1996).  Approximation of subadditive functions and convergence rates in limiting-shape results.
{\em Ann.  Probab.} {\bf 25}, 30--55.\\
Bezuidenhout, C.  and Grimmett, G.  (1991).  Exponential decay for
subcritical contact and percolation processes. {\it Ann. 
Probab.} {\bf 19}, 984--1009.\\
Chayes, J., Chayes, L. and Durrett, R. (1986).
Critical behavior of the two dimensional first passage time, {\em J.
Stat. Phys}. {\bf 45}, 933--948.\\
Cox, T. and Kesten, H. (1981). On the continuity of
the time constant of first-passage percolation, {\em J. Appl. Probab.}
{\bf 18}, 809--819.\\
Durrett, R. (1984). Oriented percolation in two
dimensions. {\em Ann. Probab.} {\bf 12}, 999--1040.\\
Durrett, R. and  Liggett, T. (1981).
The shape of the limit set in Richardson's growth model. {\em Ann.
Probab.} {\bf 9}, 186--193.\\
Grimmett, G. (1999). {\em Percolation.} Springer, Berlin.\\
 Hammersley, J. M. and  Welsh, D. J. A. (1965).
First-passage percolation, subadditive processes,
stochastic networks and generalized renewal theory.
In {\em Bernoulli, Bayes, Laplace Anniversary Volume} 
(J. Neyman   and L. LeCam,  eds.) 61--110. Springer, Berlin.\\
Kesten, H. (1986). Aspects of first-passage percolation. {\em Lecture Notes in
Math.} {\bf 1180}, 125--264. Springer, Berlin.\\
Kesten, H. and Zhang, Y. (1990). The probability of a large finite cluster in supercritical Bernoulli percolation. 
{\em Ann. Probab.} {\bf 18}, 537--555.\\
Kesten, H. and Zhang, Y. (1997). A  central limit theorem  for critical first passage percolation in two dimensions. {\em Probab. Theory  Related Fields.} {\bf 107} 137--160.\\
Marchand, R. (2002). Strict inequalities for the time
constant in first passage percolation. {\em Ann. Appl. Probab.} {\bf 12},
1001--1038.\\
Martin, J. (2004). Limit shape for directed percolation models. {\em Ann. Probab.} {\bf 32} 2908--2937.\\
Smythe, R. T. and Wierman, J. C. (1978).
First passage percolation on the square lattice.
{\em Lecture Notes in Math.} {\bf 671}. Springer, Berlin. \\
Talagrand, M. (1995). Concentration of measure and isoperimetric inequalities in product spaces. {\em Inst. Hautes Publ. Math. Etudes Sci. Publ. Math.}  {\bf 81}, 73--205.\\
Wheeden, R. and Zygmund, A. (1977). {\em Measure and integral.} Marcel Dekker, New York.\\
Yukich, J. and Zhang, Y. (2006). Singularity points for first passage percolation. {\em Ann. Probab.} {\bf 34}, 577--592. \\
Zhang, Y. (1995). Supercritical behaviors in first-passage percolation. 
{\em Stoch. Proc. Appl.} {\bf 59}, 251--266.\\
Zhang, Y. (1999). Double behavior of critical first passage percolation. In {\em Perplexing Problems in Probability}
(M. Bramson and R. Durrett eds.) 143--158. Birkhauser, Boston.\\
Zhang, Y. (2006). The divergence of fluctuations for  shape in first passage percolation. {\em Probab. Theory Related Fields.}  {\bf 136}, 298--320.\\
Zhang, Y. (2008). Shape fluctuations are different in different  directions.  {\em Ann.  Probab.} {\bf 36}, 331--362.\\

\noindent
Yu Zhang\\
Department of Mathematics\\
University of Colorado\\
Colorado Springs, CO 80933\\
email: yzhang3@uccs.edu\\
\end{document}